\documentclass[final]{doublecol}

\RequirePackage{fancyvrb}
\RequirePackage{psfrag}
\RequirePackage{algorithm}
\RequirePackage{graphicx}

\usepackage{amssymb}
\usepackage{moreverb,relsize}
\usepackage{amsbsy,natbib}
\usepackage{amsthm}
\usepackage{amsfonts,amsmath,bm}
\usepackage[dvipsone]{epsfig}
\usepackage{url}
\usepackage{natbib}
\RequirePackage{graphics,epsf}

\newcommand{\dx}{\, \mathrm{d}x}

\newcommand{\ds}{\, \mathrm{d}s}

\newcommand{\R}{\mathbb{R}}
\newcommand{\tab}{\hspace*{2em}}

\newtheorem{example}{Example}[section]

\DefineVerbatimEnvironment{code}{Verbatim}{frame=single}
\newcommand{\emp}[1]{\texttt{#1}}

\begin{document}

\LRH{Unified Form-assembly Code}

\RRH{
Martin Sandve Aln\ae{}s,
Anders Logg, Kent-Andre Mardal, Ola Skavhaug, and Hans Petter Langtangen
}

\VOL{X}

\ISSUE{X/X/X}

\PUBYEAR{2009}

\setcounter{page}{1}

\BottomCatch

\title{Unified Framework for Finite Element Assembly}

\authorA{M. S. Aln\ae{}s}
\affA{Center for Biomedical Computing, Simula Research Laboratory \\
E-mail:
martinal@simula.no\thanks{Aln\ae{}s is supported by the
Research Council of Norway under grant NFR 162730.} \\ (Corresponding author)
}

\authorB{A. Logg}
\affB{Center for Biomedical Computing, Simula Research Laboratory \\
Department of Informatics, University of Oslo \\ E-mail:
logg@simula.no\thanks{Logg is supported by an Outstanding Young Investigator grant from
the Research Council of Norway (NFR 180450).}}

\authorC{K.-A. Mardal}
\affC{Center for Biomedical Computing, Simula Research Laboratory \\
Department of Informatics, University of Oslo \\ E-mail:
kent-and@simula.no\thanks{Mardal is supported by the
Research Council of Norway under grant NFR 170650.}}

\authorD{O. Skavhaug}
\affD{Center for Biomedical Computing, Simula Research Laboratory \\
Department of Informatics, University of Oslo \\ E-mail:
skavhaug@simula.no}

\authorE{H. P.  Langtangen}
\affE{Center for Biomedical Computing, Simula Research Laboratory \\
Department of Informatics, University of Oslo \\ E-mail:
hpl@simula.no}

\begin{abstract}
Over the last fifty years, the finite element method has emerged as a
successful methodology for solving a wide range of partial
differential equations. At the heart of any finite element simulation
is the assembly of matrices and vectors from finite element
variational forms. In this paper, we present a general and unified
framework for finite element assembly. Based on this framework, we
propose a specific software interface called
\emph{Unified Form-assembly Code} (UFC)
between problem-specific and general-purpose components of finite element programs.
The interface is general in the sense that it applies to a wide range of
finite element problems (including mixed finite elements and
discontinuous Galerkin methods) and may be used with libraries that
differ widely in their design. The interface consists of a minimal set
of abstract C++ classes and data transfer is via plain C arrays.

We discuss how one may use the UFC interface to build a plug-and-play
system for finite element simulation where basic components such as
computational meshes, linear algebra and, in particular, variational
form evaluation may come from different libraries and be used
interchangeably. We further discuss how the UFC interface is used
to glue together components from the FEniCS suite of software to provide
an integrated high-level environment where variational forms may be
entered as expressions directly in Python and assembled efficiently
into sparse matrices.

A central design goal for the interface is to minimize dependency on
external libraries for the problem-specific code used in applications.
Thus, the UFC interface consists of a single C++ header file and does
not rely on external libraries for its operation. In particular, the
UFC interface does not depend on any other FEniCS components. As a
result, finite element code developers may use the interface to detach
equation-specific details from general-purpose library code, allowing
very flexible connections to alternative libraries.  We encourage
developers of finite element libraries to incorporate the interface in
their libraries. The UFC interface is released into the public domain.
\end{abstract}

\KEY{Finite elements, assembly, implementation, code generation, UFC}

\REF{to this paper should be made as follows:
M. S. Aln\ae{}s,
A. Logg, K.-A. Mardal, O. Skavhaug, and H. P. Langtangen (2008)
`Unified Framework for Finite Element Assembly'.
}

\BIO{Supported by a Center of Excellence grant from the Research
Council of Norway to the Center for Biomedical Computing at Simula Research Laboratory.}

\maketitle

\section{Introduction}

Software for solving physical problems have traditionally been
tailored to the problem at hand, often resulting in computationally
very efficient special-purpose codes. However, experience has shown that
such codes may be difficult and costly to extend to new problems.
To decrease turn-over time from problem definition to its numerical
solution, scientific code writers have to an increasingly larger
extent tried to create general libraries, containing common
numerical algorithms applicable to a wide range of problems.
Such libraries can reduce the size of the application code
dramatically and hide implementation details.
In the field of finite element solution of partial differential
equations, many general and successful
libraries have emerged during the last couple of
decades, e.g.,
Cactus,
Cogito,
COMSOL Multiphysics,
Deal.II,
Diffpack,
DOLFIN (FEniCS),
Getfem++,
Kaskade,
Sundance, and
UG (see the reference list for papers and websites).
\nocite{UG,www:Sundance,Getfempp,Cactus,Cogito96,COMSOL,www:deal.II,TCSE1,FEniCS,Kaskade,www:PETSc,Trilinos,Hypre,PyCC,www:ffc,www:dolfin}

General finite element libraries implement many standard mathematical
and numerical concepts, but the software components are often not as
carefully designed as their mathematical counterparts. From a software
engineering point of view it is important to achieve clear separation
of the various software components that build up a finite element
library, such that each component can be replaced separately. Not only
does this offer greater flexibility for application and library
developers, but it also makes the software easier to maintain,
especially under changing requirements of several developers in
long-term projects. These arguments have received much attention by
developers of general finite element libraries in recent years (see,
e.g., \cite{deal.II:paper,DUNEpaper1}).

Well designed libraries provide clear interfaces to represent this
separation.  Typically, the application code uses functions or
objects in the interface to perform basic ``high level'' steps of the
solution process. Problem-specific details, such as the variational
form, the mesh and coefficients are passed through the interface to the
library to compute a solution. Such libraries and their interfaces are
generally referred to as problem solving environments (PSEs).

However, one fundamental issue in designing such software libraries is
how to separate problem-specific code from general library
code. Some components, such as computational meshes and linear
algebra, may be implemented as reusable components (e.g. as a set of
C++ classes) with well-defined interfaces. However, other components,
such as variational forms, are intrinsically problem-specific. As a
result, those components must either be implemented and provided by
the user or \emph{generated} automatically by the
library from a high-level description of the variational form. In
either case, it becomes important to settle on a well-defined
interface for how the library should communicate with those
problem-specific components.

The design of such an interface is the subject of the present
paper. We propose a C++ interface called UFC, which provides an interface between general
reusable finite element libraries and problem-specific code.
In other libraries, the implementation of finite elements and variational
forms is usually tied to the specific mesh, matrix and vector format in use,
while in UFC we have strived to decouple these concepts.
Furthermore, the interface is designed to allow for a variety of elements
such as continuous and discontinuous Lagrange, Nedelec and Raviart-Thomas elements.

To make a successful interface, one needs a sufficiently general
framework for the underlying mathematical structures and operations.
The software interface in the current paper relies on a more general
view of variational forms and finite element assembly than commonly
found in textbooks. We therefore precisely state the mathematical
background and notation in Sections~\ref{sec:fem} and
\ref{sec:assembly}.

UFC is significantly inspired by our needs in the tools FFC, SFC, and
DOLFIN, which are software units within FEniCS,
see~\cite{KirLog2006,KirLog2007,SyFi,Log2007a}. The interplay between these
tools and UFC is explained in Section~\ref{sec:framework}, which
provides additional and more specific motivation for the design of
UFC. Highlights of the interface are covered in
Section~\ref{sec:ufc:syntax}.  Section~\ref{sec:examples} contains
some examples of high-level specifications of variational forms with
the form compilers FFC and SFC, which automatically generate code
compatible with the UFC interface for computing element matrices and
vectors. We also explain how the interface can be used with existing
libraries.

We here note that as a result of the UFC interface, the two form
compilers FFC and SFC may now be used interchangeably since both
generate code conforming to the UFC interface. These form compilers
were developed separately and independently. The work on UFC was
initially inspired by our efforts to unify the interfaces for these
form compilers.

\subsection*{Related Work}
One major reason for the success of general finite element libraries is
that many widely different physical problems can be solved by quite
short application codes utilizing the same library. The opposite
strategy, i.e., one application utilizing different alternative
libraries, has received less attention. For example, an application
might want to use an adaptive mesh data structure and its
functionality from one library, a very efficient assembly routine from
another library, basic iterative methods from, e.g.,
PETSc, combined with a preconditioner from Trilinos or Hypre.
To make this composition a true plug-and-play operation,
the various libraries would need to conform to a unified
interface to the basic operations needed in finite element solvers.
Alternatively, low level interfaces can be implemented with
thin wrapper code to connect separate software components.

In numerical linear algebra, the BLAS and LAPACK interfaces have
greatly simplified code writing. By expressing operations in the
application code in terms of BLAS and LAPACK calls, and using the
associated data (array) formats, one program can be linked to
different implementations of the BLAS and LAPACK operations. Despite
the great success of this approach, the idea has to little extent
been explored in other areas of computational science.  One recent
example is Easyviz (\cite{Easyviz}), a thin unified interface to curve
plotting and 2D/3D scalar- and vector-field visualization. This
interface allows an application program to use a MATLAB-compatible
syntax to create graphics, independently of the choice of graphics
package (Gnuplot, Grace, MATLAB, VTK, VisIt, OpenDX, etc.).  Another
example is GLAS (\cite{GLAS}), a community initiative to specify a
general interface for linear algebra libraries.  GLAS can be viewed as
an extension and modernization of the BLAS/LAPACK idea, utilizing
powerful constructs offered by C++.

Within finite elements, DUNE (\cite{DUNEpaper1,DUNEpaper2}) is
a very promising attempt to define unified interfaces between
application code and libraries for finite element computing. DUNE
provides interfaces to data structures and solution algorithms,
especially finite element meshes and iterative solution methods for
linear systems. In principle, one can write an application code
independently of the mesh data structure and the matrix solution
method.  DUNE does not directly address interfaces between the finite element
problem definition (element matrices and vectors), and the assembly
process, which is the topic of the present paper. Another difference
between DUNE and our UFC interface is the choice of programming technology
used in the interface: DUNE relies heavily on inlining via C++ templates
for efficient single-point data retrieval,
while UFC applies pointers to chunks of data. However, our view of a finite
element mesh can easily be adapted to the DUNE-Grid interface.
The DUNE-FEM module (under development) represents interfaces to various
discretization operators and serves some of the purposes of the UFC interface,
though being technically quite different.

In the finite element world, there are many competing libraries, each
with their own specialties. Thin interfaces offering only the least
common denominator functionality do not support special features for
special problems and may therefore be met with criticism. Thick
interfaces, trying to incorporate ``all'' functionality in ``all''
libraries, become too complicated to serve the original purpose of
simplifying software development.  Obtaining community consensus for
the thickness and syntax of a unified interface is obviously an
extremely challenging process. The authors of this paper suggest
another approach: a small group of people defines a thin (and hence
efficient and easy-to-use) interface, they make the software publicly
available together with a detailed documentation, and demonstrate its
advantages.  This is our aim with the present paper.

\section{Finite Element Discretization}
\label{sec:fem}

\subsection{The Finite Element}
\index{finite element}

A finite element is mathematically defined as a triplet consisting of
a polygon, a polynomial function space, and a set of linear
functionals, see~\cite{Cia78}. Given that the dimension of the
function space and the number of the (linearly independent) linear
functionals are equal, the finite element is uniquely defined. Hence,
we will refer to a finite element as a collection of
\begin{itemize}
\item a polygon $K$,
\item a polynomial space $\mathcal{P}_K$ on $K$,
\item a set of linearly independent linear functionals, the
\emph{degrees of freedom}, $\ell_i : \mathcal{P}_K \rightarrow \R, \, i =
1, 2, \ldots, n_K$.
\end{itemize}
With this definition the basis functions $\{\phi_i^K\}_{i=1}^{n_K}$ are obtained
by solving the following system of equations,
\begin{equation}
  \ell_i(\phi_j^K) = \delta_{ij}, \quad i, j = 1,2,\ldots,n_K.
\end{equation}
The computation of such a nodal basis can be automated, given (a basis
for) the polynomial space $\mathcal{P}_K$ and the set of linear
functionals $\{\ell_i\}_{i=1}^{n_K}$, see~\cite{Kir04, SyFi}.

\subsection{Variational Forms}
\index{variational form}

Consider the Poisson problem $-\nabla\cdot (w \nabla u) = f$
with Dirichlet boundary conditions on a domain $\Omega \subset
\R^d$. Multiplying by a test function $v \in V_h$ and integrating by
parts, one obtains the variational problem
\begin{equation} \label{eq:weightedpoisson}
  \int_{\Omega} w \nabla v \cdot \nabla u_h \dx = \int_{\Omega} v f \dx,
  \quad \forall v \in V_h,
\end{equation}
for the approximation $u_h \in V_h$. If $w, f \in W_h$ for some
discrete function space\footnote{It is assumed that any given function
may be represented (exactly or approximately) in some finite element
space.  Alternatively, functions may be approximated by
quadrature. Quadrature representation is not discussed here, but is
covered by the UFC specification and implemented by the form compilers
FFC and SFC.} $W_h$ we may thus write~(\ref{eq:weightedpoisson}) as
\begin{equation}
  a(v, u_h; w) = L(v; f) \quad \forall v \in V_h,
\end{equation}
where the trilinear form $a : V_h \times V_h \times W_h \rightarrow \R$ is given by
\begin{equation}
  a(v, u_h; w) = \int_{\Omega} w \nabla v \cdot \nabla u_h \dx
\end{equation}
and the bilinear form $L : V_h \times W_h \rightarrow R$ is given by
\begin{equation}
  L(v; f) = \int_{\Omega} v f \dx.
\end{equation}
Note here that $a$ is \emph{bilinear} for any given fixed $w \in W_h$
and $L$ is \emph{linear} for any given fixed $f \in W_h$.

In general, we shall be concerned with the discretization of
finite element variational forms of general arity~$r + n > 0$,
\begin{equation} \label{eq:variationalform}
  a : V_h^1 \times V_h^2 \times \cdots \times V_h^r \times
  W_h^1 \times W_h^2 \times \cdots \times W_h^n \rightarrow \R,
\end{equation}
defined on the product space $V_h^1 \times V_h^2 \times \cdots \times
V_h^r \times W_h^1 \times W_h^2 \times \cdots \times W_h^n$ of two
sets $\{V_h^j\}_{j=1}^r, \{W_h^j\}_{j=1}^n$ of discrete
function spaces on $\Omega$. We refer to
$(v_1, v_2, \ldots, v_r) \in V_h^1 \times V_h^2 \times \cdots \times V_h^r$
as \emph{primary arguments},
and to
$(w_1, w_2, \ldots, w_n) \in W_h^1 \times W_h^2 \times \cdots \times W_h^n$
as \emph{coefficients} and write
\begin{equation}
a = a(v_1, \ldots, v_r; w_1, \ldots, w_n).
\label{eq:gen_form}
\end{equation}
In the simplest case, all function spaces are equal but there are many
important examples, such as mixed methods, where the
arguments come from different function spaces.
The choice of coefficient function spaces depends on the application;
a polynomial basis simplifies exact integration, while
in some cases evaluating coefficients in quadrature points may be
required.

\subsection{Discretization}
\label{sec:Discretization}

To discretize the form $a$, we introduce a set of bases
$\{\phi_i^1\}_{i=1}^{N^1},
 \{\phi_i^2\}_{i=1}^{N^2}, \ldots,
 \{\phi_i^r\}_{i=1}^{N^r}$
for the function spaces $V_h^1, V_h^2, \ldots, V_h^r$ respectively and let $i =
(i_1, i_2, \ldots, i_r)$ be a multiindex of length $|i| = r$. The
form $a$ then defines a rank~$r$ tensor given by
\begin{equation} \label{eq:tensor}
  A_i = a(\phi_{i_1}^1, \phi_{i_2}^2, \ldots, \phi_{i_r}^r; w_1, w_2, \ldots, w_n)
  \quad \forall i \in \mathcal{I},
\end{equation}
where $\mathcal{I}$ is the index set
\begin{equation}
  \begin{split}
  & \mathcal{I} =  \prod_{j=1}^r[1,|V^j_h|] =  \\
  & \{(1,1,\ldots,1), (1,1,\ldots,2), \ldots,
  (N^1,N^2,\ldots,N^r)\}.
  \end{split}
\end{equation}
We refer to the tensor~$A$ as the \emph{discrete operator} generated
by the form~$a$ and the particular choice of basis functions.  For any
given form of arity~$r + n$, the tensor~$A$ is a (typically sparse)
tensor of rank~$r$ and dimension $|V_h^1| \times |V_h^2| \times \ldots
\times |V_h^r| = N^1 \times N^2 \times \ldots \times N^r$.
\index{global tensor}

Typically, the rank $r$ is 0, 1, or 2. When $r = 0$, the
tensor $A$ is a scalar (a tensor of rank zero), when $r = 1$, the
tensor $A$ is a vector (the ``load vector'') and when $r = 2$, the
tensor $A$ is a matrix (the ``stiffness matrix''). Forms of higher
rank also appear, though they are rarely assembled as a
higher-dimensional sparse tensor.

Note here that we consider the functions $w_1, w_2, \ldots, w_n$ as
fixed in the sense that the discrete operator~$A$ is computed for a
given set of functions, which we refer to as \emph{coefficients}. As
an example, consider again the variational
problem~(\ref{eq:weightedpoisson}) for Poisson's
equation. For the trilinear form~$a$, the rank is $r = 2$ and
the number of coefficients is $n = 1$, while for the linear form~$L$,
the rank is $r = 1$ and the number of coefficients is $n = 1$. We may
also choose to directly compute the \emph{action} of the form
$a$ obtained by assembling a vector from the form
\begin{equation}
  a(v_1; w_1, w2) = \int_{\Omega} w_1 \nabla v_1 \cdot \nabla w_2 \dx,
\end{equation}
where now $r = 1$ and $n = 2$.

We list below a few other examples to illustrate the notation.

\begin{example}
\label{example:div}
Our first example is related
to the divergence constraint in fluid flow. Let the form~$a$ be given by
\begin{equation}
a(q, u) = \int_{\Omega} q \nabla \cdot u \dx, \quad q\in V_h^1, \quad u\in V_h^2,
\end{equation}
where $V_h^1$ is a space of scalar-valued functions and where $V_h^2$
is a space of vector-valued functions.  The form $a : V_h^1 \times
V_h^2 \rightarrow \R$ has two primary arguments and thus $r = 2$.
Furthermore, the form does not depend on any coefficients and thus $n=0$.
\end{example}

\begin{example}
\label{example:linearconv}
Another common form in fluid flow (with variable density) is
\begin{equation}
a(v,u;w,\varrho) = \int_{\Omega} \varrho \, (w \cdot \nabla  u) \cdot v \dx.
\end{equation}
Here, $v\in V_h^1,\ u \in V_h^2,\ w\in W_h^1, \ \varrho \in W_h^2$, where
$V_h^1$, $V_h^2$, and $W_h^1$ are spaces of vector-valued functions, while $W_h^2$ is a space of
scalar-valued functions.
The form takes four arguments, where two of the arguments
are coefficients,
\begin{equation}
a : V_h^1 \times V_h^2 \times W_h^1 \times W_h^2 \rightarrow \R.
\end{equation}
Hence, $r=2$ and $n=2$.
\end{example}

\begin{example}
\label{example:powerlaw}
We next consider the following form appearing in nonlinear convection-diffusion with a power-law viscosity,
\begin{equation}
a(v;w,\mu,\varrho) = \int_{\Omega} \varrho (w \cdot \nabla w) \cdot v + \mu |\nabla w|^{2q} \nabla w : \nabla v \dx.
\end{equation}
Here, $v\in V_h^1,\ w\in W_h^1, \ \mu \in W_h^2, \ \varrho \in W_h^3$, where
$V_h^1$, and $W_h^1$ are spaces of vector-valued functions, while $W_h^2$ and  $W_h^3$ are spaces of scalar-valued functions.
The form takes four arguments, where three of the arguments
are coefficients,
\begin{equation}
a : V_h^1 \times W_h^1 \times W_h^2 \times W_h^3 \rightarrow \R.
\end{equation}
Hence, $r=1$ and $n=3$.
\end{example}

\begin{example}
The $H^1(\Omega)$ norm of the error $e = u - u_h$ squared is
\begin{equation}
a(;u, u_h) = \int_{\Omega} (u - u_h)^2 + |\nabla (u - u_h)|^2 \dx.
\end{equation}
The form takes two arguments and both are coefficients,
\begin{equation}
a : W_h^1 \times  W_h^2 \rightarrow \R.
\end{equation}
Hence, $r=0$ and $n=2$.
\end{example}

Defining variational forms for coupled PDEs
can be performed in two ways in the above described framework.
One approach is to couple the variational forms on the linear algebra level, using
block vectors and block matrices and defining one form for each block.
Alternatively, a single form for the coupled system may be defined using
mixed finite elements.

\section{Finite Element Assembly}
\label{sec:assembly}
\index{assembly}

The standard algorithm for computing the global sparse tensor~$A$ is
known as \emph{assembly}, see~\cite{ZieTay67,Hug87}. By this
algorithm, the tensor~$A$ may be computed by assembling (summing) the
contributions from the local entities of a finite element mesh.  To
express this algorithm for assembly of the global sparse tensor~$A$ for
a general finite element variational form of rank~$r$, we introduce
the following notation and assumptions.

Let $\mathcal{T} = \{K\}$ be a set of disjoint \emph{cells} (a
triangulation or tesselation) partitioning the domain $\Omega =
\cup_{K\in\mathcal{T}} K$. Further, let $\partial_e \mathcal{T}$
denote the set of \emph{exterior facets} (the set of cell facets
on the boundary $\partial \Omega$), and let $\partial_i
\mathcal{T}$ denote the set of $\emph{interior facets}$ (the set of
cell facets not on the boundary $\partial \Omega$).
For each discrete function space $V_h^j, \, j=1,2,\ldots,r$, we assume
that the global basis~$\{\phi_i^j\}_{i=1}^{N^j}$ is obtained by
patching together local function spaces $\mathcal{P}_K^j$ on each
cell~$K$ as determined by a local-to-global mapping.

We shall further assume that the variational
form~(\ref{eq:variationalform}) may be expressed as a sum of integrals
over the cells~$\mathcal{T}$, the exterior facets~$\partial_e
\mathcal{T}$ and the interior facets~$\partial_i \mathcal{T}$. We
shall allow integrals expressed on disjoint subsets
$\mathcal{T} = \cup_{k=1}^{n_c} \mathcal{T}_k$,
$\partial_e \mathcal{T} = \cup_{k=1}^{n_e} \partial_e \mathcal{T}_k$
and
$\partial_i \mathcal{T} = \cup_{k=1}^{n_i} \partial_i \mathcal{T}_k$
respectively.

We thus assume that the form $a$ is given by
\begin{equation}
  \begin{split}
    & a(v_1, \ldots, v_r; w_1, \ldots,  w_n) =  \\
    &\ \ \   \sum_{k=1}^{n_c} \sum_{K\in\mathcal{T}_k} \int_{K}
    I^c_k(v_1, \ldots, v_r; w_1, \ldots w_n) \dx \\
    &+
    \sum_{k=1}^{n_e} \sum_{S\in\partial_e\mathcal{T}_k} \int_{S}
    I^e_k(v_1, \ldots, v_r; w_1, \ldots,  w_n) \ds \\
    &+
    \sum_{k=1}^{n_i} \sum_{S\in\partial_i\mathcal{T}_k} \int_{S}
    I^i_k(v_1, \ldots, v_r; w_1, \ldots, w_n) \ds.
  \end{split} \label{eq:form_integrals}
\end{equation}
We refer to an integral over a cell~$K$ as a \emph{cell integral},
an integral over an exterior facet~$S$ as an \emph{exterior facet integral}
(typically used to implement Neumann and Robin type boundary conditions),
and to an integral over an interior facet~$S$ as an \emph{interior facet integral} (typically used in discontinuous Galerkin methods).

For simplicity, we consider here initially assembly of the global
sparse tensor~$A$ corresponding to a form~$a$ given by a single
integral over all cells $\mathcal{T}$, and later extend to the general
case where we must also account for contributions from several cell
integrals, interior facet integrals and exterior facet integrals.

We thus consider the form
\begin{equation}
  \begin{split}
    &a(v_1, \ldots, v_r; w_1, \ldots, w_n) = \\
    & \ \ \ \sum_{K\in\mathcal{T}} \int_K
    I^c(v_1, \ldots, v_r; w_1, \ldots, w_n) \dx,
  \end{split}
\end{equation}
for which the global sparse tensor~$A$ is given by
\begin{equation}
  A_i = \sum_{K\in\mathcal{T}} \int_K
  I^c(\phi^1_{i_1}, \ldots, \phi^r_{i_r}; w_1, \ldots, w_n) \dx.
\end{equation}
To see how to compute the tensor $A$ by summing the local
contributions from each cell~$K$, we let $n^j_K = |\mathcal{P}^j_K|$
denote the dimension of the local finite element space on $K$ for the
$j$th primary argument $v_j \in V_h^j$ for $j = 1,2,\ldots,r$. Furthermore, let
\begin{equation}
  \iota_K^j : [1,n_K^j] \rightarrow [1,N^j] \label{eq:iota_K}
\end{equation}
denote the local-to-global mapping for~$V_h^j$, that is, on any given
$K\in\mathcal{T}$, the mapping $\iota_K^j$ maps the number of a local
degree of freedom (or, equivalently, local basis function) to the
number of the corresponding global degree of freedom (or,
equivalently, global basis function). We then define for each $K \in
\mathcal{T}$ the collective local-to-global mapping $\iota_K :
\mathcal{I}_K \rightarrow \mathcal{I}$ by
\begin{equation}
  \iota_K(i) =
  (\iota_K^1(i_1),\iota_K^2(i_2),\ldots,\iota_K^r(i_r))
  \quad \forall i \in \mathcal{I}_K,
\end{equation}
where $\mathcal{I}_K$ is the index set
\begin{equation}
\begin{split}
  & \mathcal{I}_K = \prod_{j=1}^r[1,|\mathcal{P}_K^j|] \\
  & = \{(1,1,\ldots,1), (1,1,\ldots,2), \ldots,
  (n_K^1,n_K^2,\ldots,n_K^r)\}.
\end{split}
\end{equation}
Furthermore, for each $V_h^j$ we let $\{\phi^{K,j}_i\}_{i=1}^{n_K^j}$
denote the restriction to an element $K$ of the subset of the basis
$\{\phi_i^j\}_{i=1}^{N^j} \subset \mathcal{P}_K^j$ of $V_h^j$ supported on $K$.

We may now compute~$A$ by summing the contributions from
the local cells,
\begin{equation}
  \begin{split}
  A_i
  &=
  \sum_{K\in\mathcal{T}_i} \int_K
  I^c(\phi_{i_1}^1, \ldots, \phi_{i_r}^r; w_1, \ldots, w_n) \dx \\
  &=
  \sum_{K\in\mathcal{T}_i} \int_K
  I^c(\phi_{(\iota_K^1)^{-1}(i_1)}^{K,1},
      \ldots,
      \phi_{(\iota_K^r)^{-1}(i_r)}^{K,r}; w_1, \ldots, w_n) \dx \\
  &=
  \sum_{K\in\mathcal{T}_i}
  A^K_{\iota_K^{-1}(i)},
  \end{split}
\end{equation}
where $A^K$ is the local \emph{cell tensor} on cell $K$ (the ``element
stiffness matrix''), given by
\begin{equation}
  A^K_i = \int_K
  I^c(\phi_{i_1}^{K,1},
  \ldots,
  \phi_{i_r}^{K,r}; w_1, \ldots, w_n) \dx,
\end{equation}
and where $\mathcal{T}_i$ denotes the set of cells on which all basis
functions $\phi_{i_1}^1, \phi_{i_2}^2, \ldots, \phi_{i_r}^r$ are supported.
Similarly, we may sum the local contributions
from the exterior and interior facets in the form of local
\emph{exterior facet tensors} and \emph{interior facet tensors}.
\index{cell tensor}
\index{exterior facet tensor}
\index{interior facet tensor}

\begin{figure}[htbp]
  \begin{center}
    \psfrag{i0}{\hspace{-0.5cm}$\iota_K^1(1)$}
    \psfrag{i1}{\hspace{-0.5cm}$\iota_K^1(2)$}
    \psfrag{i2}{\hspace{-0.5cm}$\iota_K^1(3)$}
    \psfrag{j0}{\hspace{-0.3cm}$\iota_K^2(1)$}
    \psfrag{j1}{\hspace{-0.5cm}$\iota_K^2(2)$}
    \psfrag{j2}{\hspace{-0.1cm}$\iota_K^2(3)$}
    \psfrag{A21}{$A^K_{32}$}
    \psfrag{1}{$1$}
    \psfrag{2}{$2$}
    \psfrag{3}{$3$}
    \includegraphics[height=2in]{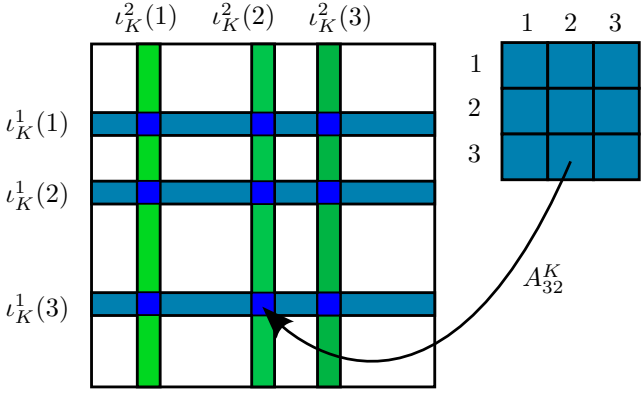}
    \caption{Adding the entries of a cell tensor~$A^K$ to the
      global tensor~$A$ using the  local-to-global mapping
      $\iota_K$, illustrated here for a rank two
      tensor (a matrix).}
    \label{fig:insertion}
  \end{center}
\end{figure}

In Algorithm~\ref{alg:assembly}, we present a general algorithm for
assembling the contributions from the local cell, exterior facet and
interior facet tensors into a global sparse tensor.  In all cases, we
iterate over all entities (cells, exterior or interior facets),
compute the local cell tensor $A^K$ (or exterior/interior facet tensor
$A^S$) and add it to the global sparse tensor as determined by the
local-to-global mapping, see~Figure~\ref{fig:insertion}.

\begin{algorithm}
$A = 0$ \\
\\
(i) \emph{Assemble contributions from all cells} \\
\textbf{for each} $K \in \mathcal{T}$ \\
\\
\tab \textbf{for} $j = 1,2,\ldots,r$: \\
\tab\tab Tabulate the local-to-global mapping $\iota_K^j$ \\
\\
\tab \textbf{for} $j = 1,2,\ldots,n$: \\
\tab\tab Extract the values of $w_j$ on $K$
\\
\\
\tab Take $0 \leq k \leq n_c$ such that $K \in \mathcal{T}_k$ \\
\tab Tabulate the cell tensor $A^K$ for $I^c_k$ \\
\tab Add $A^K_i$ to $A_{\iota_K^1(i_1), \iota_K^2(i_2), \ldots, \iota_K^r(i_r)}$ for $i\in I_K$ \\
\\
(ii) \emph{Assemble contributions from all exterior facets} \\
\textbf{for each} $S \in \partial_e\mathcal{T}$ \\
\\
\tab \textbf{for} $j = 1,2,\ldots,r$: \\
\tab\tab Tabulate the local-to-global mapping $\iota_{K(S)}^j$ \\
\\
\tab \textbf{for} $j = 1,2,\ldots,n$: \\
\tab\tab Extract the values of $w_j$ on $K(S)$
\\
\\
\tab Take $0 \leq k \leq n_e$ such that $S \in \partial_e \mathcal{T}_k$ \\
\tab Tabulate the exterior facet tensor $A^S$ for $I^e_k$ \\
\tab Add $A^S_i$ to $A_{\iota_{K(S)}^1(i_1), \iota_{K(S)}^2(i_2), \ldots, \iota_{K(S)}^r(i_r)}$ for $i\in I_{K(S)}$ \\
\\
\\
(iii) \emph{Assemble contributions from all interior facets} \\
\textbf{for each} $S \in \partial_i\mathcal{T}$ \\
\\
\tab \textbf{for} $j = 1,2,\ldots,r$: \\
\tab\tab Tabulate the local-to-global mapping $\iota_{K(S)}^j$ \\
\\
\tab \textbf{for} $j = 1,2,\ldots,n$: \\
\tab\tab Extract the values of $w_j$ on $K(S)$
\\
\\
\tab Take $0 \leq k \leq n_i$ such that $S \in \partial_i \mathcal{T}_k$ \\
\tab Tabulate the interior facet tensor $A^S$ for $I^i_k$ \\
\tab Add $A^S_i$ to $A_{\iota_{K(S)}^1(i_1), \iota_{K(S)}^2(i_2), \ldots, \iota_{K(S)}^r(i_r)}$ for $i\in I_{K(S)}$ \\
\caption{Assembling the global tensor~$A$ from the local contributions
  on all cells, exterior and interior facets. For assembly over
  exterior facets, $K(S)$ refers to the cell $K\in\mathcal{T}$ incident
  with the exterior facet~$S$, and for assembly over interior facets,
  $K(S)$ refers to the ``macro cell'' consisting of the pair of cells
  $K^+$ and $K^-$ incident with the interior facet~$S$.}
\label{alg:assembly}
\end{algorithm}

\section{Software Framework for Finite Element Assembly}
\label{sec:framework}

In a finite element application code, typical input from the user is the
variational (weak) form of a PDE, a choice of finite elements, a
geometry represented by a mesh, and user-defined functions that appear
as coefficients in the variational form. For a linear PDE, the typical
solution procedure consists of first assembling a (sparse) linear
system $AU = b$ from given user input and then solving that linear
system to obtain the degrees of freedom~$U$ for the discrete finite
element approximation~$u_h$ of the exact solution~$u$ of the PDE.
Even when the solution procedure is more involved, as for a nonlinear
problem requiring an iterative procedure, each iteration may involve
assembling matrices and vectors. It is therefore clear that the
assembly of matrices and vectors (or in general tensors) is an
important task for any finite element software framework. We refer to
the software component responsible for assembling a global tensor
from given user input consisting of a variational form, finite element
function spaces, mesh and coefficients as the~\emph{Assembler}.

As demonstrated in Algorithm~\ref{alg:assembly}, the Assembler needs
to iterate over the cells in the mesh, tabulate degree of freedom
mappings, extract local values of coefficients, compute the local
element tensor, and add each element tensor to the global tensor which
is the final output. Thus, the Assembler is a software component where
many other components are combined. It is therefore important that the
software components on which the Assembler depends have well-defined
interfaces. We discuss some issues relating to the design of these
software components below and then demonstrate how these software
components together with the Assembler may be combined into a general software
framework for finite element assembly.

\subsection{Variational Forms}

Implementations of discrete variational forms in a general
finite element library usually consist of programming expressions for
the integrands $I^c_k$, $I^e_k$, $I^i_k$ (see
Equation~\eqref{eq:form_integrals}), eventually writing a quadrature
loop and a loop over element matrix indices depending on the abstraction
level of the library (see \cite{deal.II:paper,SC.1.Langtangen.2003}).
An alternative approach is to apply exact integration instead of
quadrature. In either case, the result of this computation may be
communicated through the UFC interface.

The motivation behind the UFC interface is to separate the implementation
of the form from other details of the assembly such as the mesh and the
linear algebra libraries in use.

In the FEniCS finite element software framework,
a high-level form language embedded in Python is used to define
the variational form and finite elements.
This reduces the distance from the mathematical formulation of a PDE to
an implementation of a PDE solver, removes tedious and error-prone tasks
such as coding PDE-specific assembly loops, and enables rapid
prototyping of new models and methods.
To retain computational efficiency, we generate efficient low-level
code from the abstract form description, using exact integration
where possible. Code generation adds another
complexity layer to the software, and it becomes even more important
to keep a clear separation between software components such that the
interface between generated code and library code is well defined.
This is achieved by generating implementations of the UFC interface.

\subsection{Mesh Libraries}

Many different representations of computational meshes exist.
Typically, each finite element library provides its own internal
implementation of a computational mesh. We do not wish to tie the
UFC interface to one particular mesh representation or one particular
library. Still, several operations like the element tensor computation
depends on local mesh data. For this reason, the UFC interface
provides a low-level data structure to communicate single cell
data. In addition, a small data structure is used to communicate
global mesh dimensions which are necessary for computing the
local-to-global mapping. Assemblers implemented on top of the UFC
interface must therefore be able to copy/translate cell data from the
mesh library being used to the UFC data structure (involving a
minimal overhead). This makes it possible to achieve separation
between the mesh representation and the element tensor
computation. The Assembler implemented in FEniCS (as part of DOLFIN)
is implemented for one particular mesh format,
see~\cite{submitted_Log2008}, but an Assembler component could easily
be written for other mesh libraries like the PETSc~\emph{Sieve},
see~\cite{KneKar05} and the DUNE-Grid interface
(\cite{DUNEpaper1,DUNEpaper2}).

\subsection{Linear Algebra Libraries}

It is desirable to reuse existing high-performance linear algebra
libraries like PETSc and Trilinos. It is therefore important that the
Assembler is able to assemble element tensors into global matrices and
vectors implemented by external libraries. Aggregation into matrix and
vector data structures are fairly similar operations, typically
consisting of passing array pointers to existing functions in the
linear algebra libraries.  By implementing a common interface for
assembly into tensors of arbitrary rank, the same assembly routine can
be reused for any linear algebra library without changes. This avoids
duplication of assembly code, and one may easily change the output
format of the assembly procedure. The details of these interfaces are
beyond the scope of the current paper.  At the time of writing, we
have written assembly routines with support for matrices and vectors
from Epetra (Trilinos), PETSc, PyCC, and uBLAS in addition to scalars.

\begin{figure}[H]
  \begin{center}
    \includegraphics[width=0.45\textwidth]{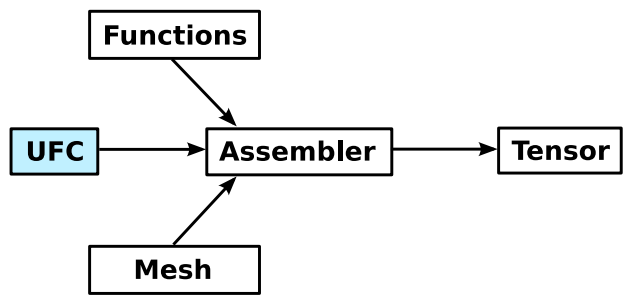}
    \caption{Assembling a tensor from a given UFC, mesh and functions (coefficients).}
    \label{fig:inoutassembler}
  \end{center}
\end{figure}

With components available for finite element variational forms, mesh
representation, and linear algebra, we may use the UFC interface to
combine these components to build a PSE for partial differential
equations. The central component of this PSE is the Assembler. As
illustrated in Figure~\ref{fig:inoutassembler}, the Assembler takes as input a
variational form, communicated through the UFC interface, a mesh and
a set of functions (the coefficients), and assembles a tensor.

\subsection{High-level Interfaces}

In FEniCS, we have additional application-level abstractions for
expressing variational forms, meshes, functions and linear algebra
objects to achieve a consistent high-level user-interface.
The generation of the UFC code may then be hidden from the user,
who just provides a high-level description of the form. The PSE may
then automatically generate the UFC at run-time, functioning as
just-in-time (JIT) compiler, and call the Assembler with the generated
UFC. Below, we demonstrate how this may be done in the Python
interface of DOLFIN. The user here defines a finite element function
space, and a pair of bilinear and linear forms $a(v, u) =
\int_{\Omega} \nabla v \cdot \nabla u + vu \dx$ and $L = \int_{\Omega}
vf \dx$, from which a matrix and vector may be assembled by calls to
the function~\texttt{assemble}. A linear system solver may then be invoked
to compute the degrees of freedom $U$ of the solution.

{\scriptsize
\begin{code}
element = FiniteElement("CG", "triangle", 1)

v = TestFunction(element)
u = TrialFunction(element)
f = Function(element, mesh, 100.0)

a = dot(grad(v), grad(u))*dx + v*u*dx
L = v*f*dx

A = assemble(a, mesh)
b = assemble(L, mesh)

U = solve(A, b)
\end{code}
}

\begin{figure}[H]
  \begin{center}
  \includegraphics[width=0.47\textwidth]{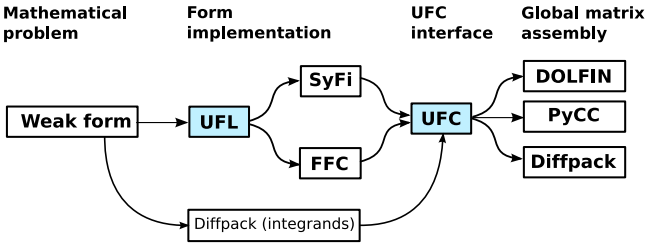}
  \caption{Alternate routes from mathematical description to matrix
    assembly enabled by the UFC interface (Note that the Diffpack example is fictional).}
    \label{fig:softwarealternatives}
  \end{center}
\end{figure}

While FEniCS provides an integrated environment, including the PSE
DOLFIN and the two form compilers FFC (FEniCS Form Compiler) and SFC (SyFi Form Compiler),
the fact that these components comply with the UFC interface means that they may also be
used interchangeably in heterogeneous environments together with other
libraries (that implement or use the UFC interface). This is illustrated in
the flow diagram of Figure~\ref{fig:softwarealternatives} where
alternate routes from mathematical description to matrix assembly are
demonstrated.

\begin{figure}[H]
  \begin{center}
    \includegraphics[width=0.45\textwidth]{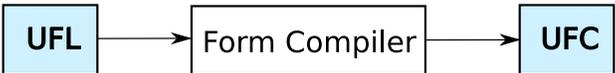}
    \caption{An abstract definition (UFL) of a finite element
      variational form is given as input to a form compiler, which
      generates UFC code as output.}
    \label{fig:uflufc}
  \end{center}
\end{figure}

In Figure~\ref{fig:softwarealternatives}, we have also included
another interface UFL (Unified Form Language) which provides a unified
way to express finite element variational forms.  The UFL interface is
currently in development. Together, UFL and UFC provide a unified
interface for the input and output of form compilers, see
Figure~\ref{fig:uflufc}
.

\section{The UFC Interface}
\label{sec:ufc:syntax}

\begin{figure}[H]
\includegraphics[width=0.45\textwidth]{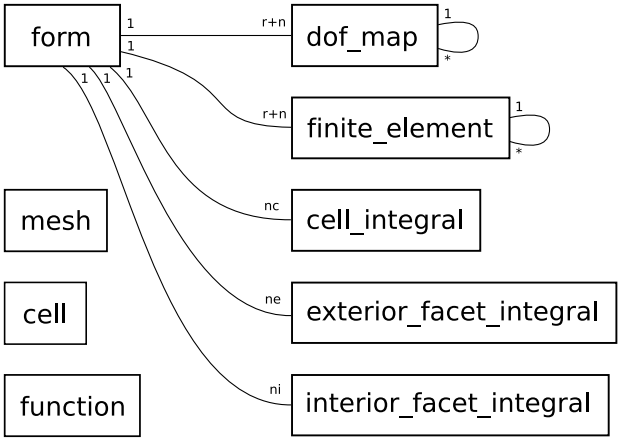}
\caption{UML diagram of the UFC class relations}
\label{fig:uml}
\end{figure}

The UFC interface consists of a small collection of abstract C++
classes that represent common components for assembling tensors using
the finite element method. These classes are accompanied by a well
documented (\cite{AlnLan2008}) set of conventions for numbering of
cell data and other arrays. We have strived to make the classes as
simple as possible while not sacrificing generality or efficiency.
Data is passed as plain C arrays for efficiency and minimal
dependencies.  Most functions are declared \emp{const}, reflecting
that the operations they represent should not change the outcome of
future operations.\footnote{The exceptions are the functions to
initialize a \emp{dof\_map}.}  Other initialization of
implementation-specific data should ideally be performed in
constructors.

One can ask why the UFC interface consists of classes and not
plain functions. There are three reasons for this.
First, we want to handle each form as a self contained ``black box'',
which can be passed around easily in an application.
Many functions belong together conceptually, thus making it
natural to collect them in a class ``namespace''.
Second, we need multiple versions of each function in the software
representation of variational forms, in particular to represent
multiple variational forms and multiple finite element function spaces.
This is best achieved by making each such function a member function
of a class and having multiple instances of that class.
Third, UFC function implementations may need access to stored data,
and with a plain function-based interface these data would then
need to be global variables.
In particular, when existing libraries or applications want to implement the
UFC interface, it may be necessary for the subclasses of UFC classes
to inherit from existing classes or to have pointers to other objects.

\subsection{Class Relations}

Figure (\ref{fig:uml}) shows all the classes and their relations.  The
classes \emp{mesh}, \emp{cell}, and \emp{function} provide the means
for communicating mesh and coefficient function data as
arguments. Each argument of the form (both primary arguments and
coefficients) is represented by a \emp{finite\_element} and
\emp{dof\_map} object. The integrals are represented by one of the
classes \emp{cell\_integral}, \emp{exterior\_facet\_integral}, or
\emp{interior\_facet\_integral}. An object of the class \emp{form}
gives access to all other objects in a particular implementation.
In this paper, we will not describe all the functions of these classes
in detail. A complete specification can be found in the
manual~(\cite{AlnLan2008}).

At the core of UFC is the class \emp{form}, which represents the general
variational form~$a$ of Equation~\eqref{eq:form_integrals}.
Subclasses of \emp{form} must implement factory
functions which may be called to create \emp{cell\_integral},
\emp{exterior\_facet\_integral} and
\emp{interior\_facet\_integral} objects. These objects in turn
know how to compute their respective contribution from a cell or facet
during assembly. A code fragment from the \emp{form} class
declaration is shown below.

{\scriptsize
\begin{code}
class form
{
public:

  ...

  /// Create cell integral on sub domain i
  virtual cell_integral*
  create_cell_integral(unsigned int i)
  const = 0;

  /// Create exterior facet integral on sub domain i
  virtual exterior_facet_integral*
  create_exterior_facet_integral(unsigned int i)
  const = 0;

  /// Create interior facet integral on sub domain i
  virtual interior_facet_integral*
  create_interior_facet_integral(unsigned int i)
  const = 0;
};
\end{code}
}

The \emp{form} class also specifies functions for creating
\emp{finite\_element} and \emp{dof\_map} objects for the finite
element function spaces $\{V_h^j\}_{j=1}^r, \{W_h^j\}_{j=1}^n$ of the
variational form.
The \emp{finite\_element} object provides functionality
such as evaluation of degrees of freedom
and evaluation of basis functions and their derivatives.
The \emp{dof\_map} object provides functionality such as tabulating the
local-to-global mapping of degrees of freedom on a single element,
as well as tabulation of subsets associated with particular mesh entities, used
to apply Dirichlet boundary conditions and build connectivity information.

Both the \emp{finite\_element} and \emp{dof\_map} classes can
represent mixed elements, in which case it is possible to obtain
\emp{finite\_element} and \emp{dof\_map} objects for each sub-element
in a hierarchical manner.  Vector elements composed of scalar elements
are in this context seen as special cases of mixed elements where
all sub-elements are equal.  Thus, e.g., from a \emp{dof\_map}
representing a $P_2-P_1$ Taylor-Hood element, it is possible to
extract one \emp{dof\_map} for the quadratic vector element and one
\emp{dof\_map} for the linear scalar element. From the vector element,
a \emp{dof\_map} for the quadratic scalar element of each vector
component can be obtained. This can be used to access subcomponents
from the solution of a mixed system.

\subsection{Stages in the Assembly Algorithm}

\begin{figure}[H]
{\scriptsize
\begin{code}
enum shape {interval, triangle, quadrilateral,
            tetrahedron, hexahedron};

class cell {
public:
  shape cell_shape;
  unsigned int topological_dimension;
  unsigned int geometric_dimension;

  /// Array of global indices for the mesh entities of the cell
  unsigned int** entity_indices;

  /// Array of coordinates for the vertices of the cell
  double** coordinates;
};
\end{code}
}
\caption{Data structure for communicating single cell data.}
\label{fig:cellcode}
\end{figure}

Next, we focus on a few key parts of
the interface and explain how these can be used to implement the
assembly algorithm (Algorithm \ref{alg:assembly}). This algorithm
consists of three stages: (i) assembling the contributions from all
cells, (ii) assembling the contributions from all exterior facets, and
(iii) assembling the contributions from all interior facets.

Each of the three assembly stages (i)--(iii) of
Algorithm~\ref{alg:assembly} is further composed of five steps. In the
first step, the polygon $K$ is fetched from the mesh, typically
implemented by filling a \emp{cell} structure (see Figure \ref{fig:cellcode}) with coordinate
data and global numbering of the mesh entities in the cell.
This step depends on the specific mesh being used.

Secondly, the local-to-global mapping of degrees of freedom is tabulated for
each of the function spaces.  That is, for each of the discrete finite
element spaces on $K$, we tabulate (or
possibly compute) the global indices for the degrees of freedom
on $\{V_h^j\}_{j=1}^r$ and $\{W_h^j\}_{j=1}^n$.

The class \emp{dof\_map} represents the mapping between local and
global degrees of freedom for a finite element space. A
\emp{dof\_map} is initialized with global mesh dimensions by
calling the function \emp{init\_mesh(const mesh\& m)}. If this
function returns \emp{true}, the \emp{dof\_map} should be additionally
initialized by calling the function
\emp{init\_cell(const mesh\& m, const cell\& c)}
for each cell in the global mesh, followed by
\emp{init\_cell\_finalize} after the last cell.
After the initialization stage, the
mapping may be tabulated at a given cell by calling a function with
the following signature.

{\scriptsize
\begin{code}
void dof_map::tabulate_dofs(unsigned int* dofs,
                            const mesh& m,
                            const cell& c) const
\end{code}
}
Here, \emp{unsigned int* dofs} is a pointer to the first element of an
array of unsigned integers that will be filled with the
local-to-global mapping on the current cell during the function call.

In the third step of each stage of Algorithm~\ref{alg:assembly}, we
may use the tabulated local-to-global mapping to interpolate (extract)
the local values of any of the coefficients in
$\{W_h^j\}_{j=1}^n$.

If a coefficient $w_j$ is not given as a linear combination of
basis functions for $W_h^j$, it must at this
step be interpolated into $W_h^j$, using the interpolant defined by
the degrees of freedom of $W_h^j$ (for example point evaluation at a
set of nodal points). In this case, the coefficient function is passed
as an implementation of the \emp{function} interface (a simple functor)
to the function~\texttt{evaluate\_dofs}.
{\scriptsize
\begin{code}
/// Evaluate linear functionals for all dofs on the function f
void finite_element::evaluate_dofs(double *values,
                                   const function& f,
                                   const cell& c) const
\end{code}
}

In the fourth step, the local element tensor contributions (cell or
exterior/interior facet tensors) are computed. This is done by a call
to the function~\texttt{tabulate\_tensor}, illustrated below for a cell
integral.

\scriptsize
\begin{code}
void cell_integral::tabulate_tensor(double* A,
                                    const double * const * w,
                                    const cell& c) const
\end{code}
\normalsize
Similarly, one may evaluate interior and exterior facet contributions
using slightly different function signatures.

Finally, at the fifth step, the local element tensor contributions are added to the global tensor,
using the local-to-global mappings previously obtained by calls to the
\emp{tabulate\_dofs} function. This is a simple operation that depends on the linear algebra library in use.

\section{Examples}
\label{sec:examples}

In this section, we demonstrate how UFC is used in practice in
DOLFIN, FFC, and SFC. First, we show a part of the assembly algorithm
(Algorithm~\ref{alg:assembly}) as implemented in DOLFIN. We then show
examples of input to the form compilers FFC and SFC as well as part of
the corresponding UFC code generated as output.
Examples include Poisson's equation and linear convection (see
Example \ref{example:linearconv}).

\subsection{An Example UFC Assembler}

To demonstrate how one may implement an assembler based on the UFC
interface, we provide here a (somewhat simplified) excerpt from the
DOLFIN assembler.\footnote{The \texttt{ufc} object is here an instance
of a simple DOLFIN class that holds pointers to arrays and UFC
container classes, such as the array \texttt{A} and cell data
\texttt{ufc::cell}, needed to communicate through the UFC interface.}

{\scriptsize
\begin{code}
for (CellIterator cell(mesh); !cell.end(); ++cell)
{
  ufc.update(*cell);

  for (uint i = 0; i < ufc.form.rank(); i++)
    dof_map_set[i].tabulate_dofs(ufc.dofs[i], *cell);

  for (uint i = 0; i < coefficients.size(); i++)
    coefficients[i]->interpolate(ufc.w[i], ufc.cell,
                                 *ufc.coefficient_elements[i],
                                 *cell);

  integral->tabulate_tensor(ufc.A, ufc.w, ufc.cell);

  A.add(ufc.A, ufc.local_dimensions, ufc.dofs);
}
\end{code}
}

The outer loop iterates over all cells in a given mesh. For each cell,
a \texttt{ufc::cell} is updated and the local-to-global mapping is constructed.
We then interpolate all the form coefficients on the cell and compute the
element tensor. At the end of the iteration, the local-to-global mapping is used
to add the local tensor to the global tensor.

\subsection{FFC Examples}
The form compiler FFC provides a simple language for specification of
variational forms, which may be entered either directly in Python or
in text files given to the compiler on the command-line. For each
variational form given as input, FFC generates UFC-compliant C++ code
for evaluation of the corresponding element tensor(s).

\subsubsection*{Poisson's Equation}
We consider the following input file to FFC for Poisson's equation.
{\scriptsize
\begin{code}
element = FiniteElement("Lagrange", "triangle", 1)

v = TestFunction(element)
u = TrialFunction(element)
f = Function(element)

a = dot(grad(v), grad(u))*dx
L = v*f*dx
\end{code}
}

Here, two forms \texttt{a} (bilinear) and \texttt{L} (linear) are defined.
Both the test and trial spaces are spanned by linear Lagrange elements
on triangles in two dimensions. When compiling this code using FFC, a C++
header file is created, containing UFC code that may be used to
assemble the global sparse stiffness matrix and load vector. Below, we
present the code generated for evaluation of the element stiffness
matrix for the bilinear form~\texttt{a}.

{\scriptsize
\begin{code}
virtual void tabulate_tensor(double* A,
                             const double * const * w,
                             const ufc::cell& c) const
{
  // Extract vertex coordinates
  const double * const * x = c.coordinates;

  // Compute Jacobian of affine map from reference cell
  const double J_00 = x[1][0] - x[0][0];
  const double J_01 = x[2][0] - x[0][0];
  const double J_10 = x[1][1] - x[0][1];
  const double J_11 = x[2][1] - x[0][1];

  // Compute determinant of Jacobian
  double detJ = J_00*J_11 - J_01*J_10;

  // Compute inverse of Jacobian
  const double Jinv_00 =  J_11 / detJ;
  const double Jinv_01 = -J_01 / detJ;
  const double Jinv_10 = -J_10 / detJ;
  const double Jinv_11 =  J_00 / detJ;

  // Set scale factor
  const double det = std::abs(detJ);

  // Compute geometry tensors
  const double G0_0_0 = det*(Jinv_00*Jinv_00 + Jinv_01*Jinv_01);
  const double G0_0_1 = det*(Jinv_00*Jinv_10 + Jinv_01*Jinv_11);
  const double G0_1_0 = det*(Jinv_10*Jinv_00 + Jinv_11*Jinv_01);
  const double G0_1_1 = det*(Jinv_10*Jinv_10 + Jinv_11*Jinv_11);

  // Compute element tensor
  A[0] = 0.5*G0_0_0 + 0.5*G0_0_1 + 0.5*G0_1_0 + 0.5*G0_1_1;
  A[1] = -0.5*G0_0_0 - 0.5*G0_1_0;
  A[2] = -0.5*G0_0_1 - 0.5*G0_1_1;
  A[3] = -0.5*G0_0_0 - 0.5*G0_0_1;
  A[4] = 0.5*G0_0_0;
  A[5] = 0.5*G0_0_1;
  A[6] = -0.5*G0_1_0 - 0.5*G0_1_1;
  A[7] = 0.5*G0_1_0;
  A[8] = 0.5*G0_1_1;
}
\end{code}
}
In FFC, an element tensor contribution is computed as a tensor
contraction between a geometry tensor varying from cell to cell, and a
geometry independent tensor on a reference element, see
\cite{KirLog2006,KirLog2007}. For simple forms, like the one under
discussion, the main work is then to construct the geometry tensor,
related to the geometrical mapping between the reference element and
physical element.

Having computed the element tensor, one needs to compute the
local-to-global mapping in order to know where to insert the local
contributions in the global tensor. This mapping may be obtained by
calling the member function \texttt{tabulate\_dofs} of the class
\texttt{ufc::dof\_map}. FFC uses an implicit ordering scheme, based on
the indices of the topological entities in the mesh. This information
may be extracted from the \texttt{ufc::cell} attribute
\texttt{entity\_indices}.

{\scriptsize
\begin{code}
virtual void tabulate_dofs(unsigned int* dofs,
                           const ufc::mesh& m,
                           const ufc::cell& c) const
{
  dofs[0] = c.entity_indices[0][0];
  dofs[1] = c.entity_indices[0][1];
  dofs[2] = c.entity_indices[0][2];
}
\end{code}
}

For Lagrange elements on triangles, each degree of freedom is
associated with a global vertex. Hence, FFC constructs the mapping by
picking the corresponding global vertex number for each degree of
freedom.

\subsubsection*{Linear Convection}

Consider the variational form in Example \ref{example:linearconv}. The input
file to FFC reads as follows.

{\scriptsize
\begin{code}
vector_element = VectorElement("Lagrange", "triangle", 1)
scalar_element = FiniteElement("Lagrange", "triangle", 1)

v = TestFunction(vector_element)
u = TrialFunction(vector_element)
w = Function(vector_element)
rho = Function(scalar_element)

a = rho*v[i]*w[j]*u[i].dx(j)*dx
\end{code}
}

The code generated for the \texttt{tabulate\_tensor} function is
presented below. Computations involving coefficients are performed by
interpolating the functions $w$ and $\rho$ on the cell under
consideration. These values are stored in the array \texttt{w}
below. For clarity, some code has been omitted in this example.

{\scriptsize
\begin{code}
virtual void tabulate_tensor(double* A,
                             const double * const * w,
                             const ufc::cell& c) const
{
  // Extract vertex coordinates and compute Jacobian etc
  // as in previous example
  ...
  // Compute coefficients
  const double c1_0_0_0 = w[1][0];
  const double c1_0_0_1 = w[1][1];
  const double c1_0_0_2 = w[1][2];
  const double c0_0_1_0 = w[0][0];
  ...
  const double c0_0_1_5 = w[0][5];
  // Compute geometry tensors
  const double G0_0_0_0_0 = det*c1_0_0_0*c0_0_1_0*Jinv_00;
  const double G0_0_0_1_0 = det*c1_0_0_0*c0_0_1_0*Jinv_10;
  const double G0_0_1_0_0 = det*c1_0_0_0*c0_0_1_1*Jinv_00;
  ...
  // Compute element tensor
  A[0] = -0.05*G0_0_0_0_0 - ...
  A[1] =  0.05*G0_0_0_0_0 + ...
  A[2] =  0.05*G0_0_0_1_0 + ...
  ...
}
\end{code}
}

The local-to-global mapping for the space of piecewise linear vectors
is computed by associating two values with each vertex. The code
generated for \texttt{tabulate\_dofs} is presented below.

{\scriptsize
\begin{code}
virtual void tabulate_dofs(unsigned int* dofs,
                           const ufc::mesh& m,
                           const ufc::cell& c) const
{
  dofs[0] = c.entity_indices[0][0];
  dofs[1] = c.entity_indices[0][1];
  dofs[2] = c.entity_indices[0][2];
  unsigned int offset = m.num_entities[0];
  dofs[3] = offset + c.entity_indices[0][0];
  dofs[4] = offset + c.entity_indices[0][1];
  dofs[5] = offset + c.entity_indices[0][2];
}
\end{code}
}

FFC generates code for arbitrary multilinear forms and currently
supports arbitrary degree continuous Lagrange elements, discontinuous
elements, RT elements, BDM elements, BDFM elements and Nedelec
elements in two and three space dimensions.

\subsection{SFC Examples}

SFC is another form compiler producing UFC code, in which the user
defines variational forms in Python using a symbolic engine based on
GiNaC~(\cite{www:ginac}).  It has a slightly different feature set
than FFC, such as using symbolic differentiation to automatically
compute the Jacobi matrix of a nonlinear form. The resulting low-level
UFC code is very similar.

\subsubsection*{Power-law Viscosity}

Example \ref{example:powerlaw} is specified in SFC as follows.

{\scriptsize
\begin{code}
vector_element = VectorElement("Lagrange", "triangle", 1)
scalar_element = FiniteElement("DG", "triangle", 0)

v = TestFunction(vector_element)
w = Function(vector_element)
mu = Function(scalar_element)
rho = Function(scalar_element)

def power_law(v, w, mu, rho, itg):
    q = 0.3
    GinvT = itg.GinvT()
    Dw  = grad(w, GinvT)
    Dv  = grad(v, GinvT)
    wDw = dot(w, Dw)

    return rho*dot(wDw, v) + mu*inner(Dw,Dw)**q * inner(Dw,Dv)

F_form = Form(basisfunctions = [v],
              coefficients   = [w, mu, rho])
F_form.add_cell_integral(power_law)
J_form = Jacobi(F_form)
\end{code}
}

The syntax for defining elements and arguments is the same
as in FFC, but the integrand is specified in a slightly different syntax\footnote{The variable \texttt{itg} is an integral object containing information about the mapping between physical
coordinates and the reference element.}.
This code also computes the form corresponding to the
Jacobian matrix, using symbolic differentiation.
The generated code for computing the Jacobian matrix is in this
case more complicated but it implements
the UFC interface in the same manner as in the previous examples.

In the current implementation, SFC explicitly constructs a
local-to-global mapping at run-time.
In this case, with Lagrange elements,  the global coordinates
identify the degrees of freedom. The UFC interface supports
constructing the local-to-global mapping  through
the \texttt{init\_mesh} and \texttt{init\_cell} methods of
\texttt{ufc::dof\_map}. Below, we present the code generated for
\texttt{init\_cell} (where we use
additional structures of type \emp{Ptv} (point) for representing degrees of freedom  and the container \emp{Dof\_Ptv dof} for building the  local-to-global mapping).

{\scriptsize
\begin{code}
void dof_map_2D::init_cell(const ufc::mesh& m, const ufc::cell& c)
{
  // coordinates
  double x0 = c.coordinates[0][0]; double y0 = c.coordinates[0][1];
  double x1 = c.coordinates[1][0]; double y1 = c.coordinates[1][1];
  double x2 = c.coordinates[2][0]; double y2 = c.coordinates[2][1];

  // affine map
  double G00 = x1 - x0;
  double G01 = x2 - x0;

  double G10 = y1 - y0;
  double G11 = y2 - y0;

  unsigned int element = c.entity_indices[2][0];

  double dof0[2] = { x0, y0 };
  Ptv pdof0(2, dof0);
  dof.insert_dof(element, 0, pdof0);

  double dof1[2] = {  G00+x0,  y0+G10 };
  Ptv pdof1(2, dof1);
  dof.insert_dof(element, 1, pdof1);

  double dof2[2] = {  x0+G01,  G11+y0 };
  Ptv pdof2(2, dof2);
  dof.insert_dof(element, 2, pdof2);
}
\end{code}
}

The \texttt{dof\_map} class is only responsible for the uniqueness of
the local-to-global mapping. Possible renumbering strategies may be
imposed by the assembler, for example to minimize communication when
assembling in parallel.

\section{Discussion}

We have used (generated) UFC for many applications, including
Poisson's equation; convection--diffusion--reaction equations;
continuum equations for linear elasticity, hyperelasticity, and
plasticity; the incompressible Navier-Stokes equations; and mixed
formulations for the Hodge Laplacian. The types of finite elements
involved include standard continuous Lagrange elements of arbitrary
order, discontinuous Galerkin formulations, BDM elements,
Raviart--Thomas elements, Crouzeix--Raviart elements, and Nedelec
elements.

The form compilers FFC and SFC are UFC compliant, both generating
efficient UFC code from an abstract problem definition. Assemblers
have been implemented in DOLFIN and PyCC, using the DOLFIN
mesh representation, and together covering linear algebra formats from
PETSc, Trilinos (Epetra), uBLAS, and PyCC.  Parallel assembly is
supported in the current development version of DOLFIN, without
requiring any modifications to UFC since it operates on an element
level.  Altogether, this demonstrates that the UFC interface is
flexible both in terms of the applications and finite element
formulations it covers, and in terms of its interoperability with
existing libraries.

One of the main limitations in the current version of the UFC
interface (v1.1) is the assumption of a homogeneous mesh, that is,
only one cell shape is allowed throughout the mesh. Thus, although
mesh ordering conventions have been defined for the interval,
triangle, tetrahedron, quadrilateral, and hexahedron, only one type of
shape can be used at any time.  Also, higher order (non-affinely
mapped) meshes are not supported in the current version of the
interface.  Another limitation is that only one fixed finite element
space can be chosen for each argument of the form, which excludes
$p$-refinement (increasing the element order in a subset of the
cells).  All these limitations may be removed in future versions of
UFC, and we encourage interested developers to make contact to address
these limitations.

UFC provides a unified interface for code generated as \emph{output}
by form compilers such as FFC and SFC. Similarly, we are currently
working on a specification for a unified form language (UFL) to
function as a common \emph{input} to form compilers. Currently, both
FFC and SFC provide (different) form languages for easy specification
of variational forms in a high-level syntax. With a unified form
language, a user may specify a variational form in that language and
assemble the corresponding discrete operator (tensor), independently
of the components being used to generate the UFC code from the UFL, and
independently of the components being used to assemble the tensor from
the UFC form.

\section{Conclusion}

We have presented a general framework for assembly of finite element
variational forms. Based on this framework, we have then extracted an
interface (UFC) that may be used to provide a communication layer
between general-purpose and problem-specific code for assembly of
finite element variational forms.

The interface makes minimal assumptions on the type of problem being
solved and the data structures involved. For example, the discrete
variational form may in general be multilinear and hence assemble into
a tensor of arbitrary rank. The basic data structures used to pass
data through the interface are composed of plain C arrays.  The
minimal set of assumptions on problem and data structures enables
application of the interface to a wide range of variational forms and
a large collection of finite element libraries.

We have used the UFC interface in the implementation of the FEniCS
suite of finite element tools. In a simple Python script, one may
define a variational form and a mesh, and assemble the corresponding
global sparse matrix (or vector). When doing so, UFC code is generated
by either of the form compilers FFC or SFC, and passed to the
UFC-compatible assembler of the general-purpose finite element library
DOLFIN.

We encourage developers of finite element software to use the UFC
interface in their libraries. By doing so, those libraries may
directly take advantage of the form compilers FFC or SFC to specify
finite element problems. Moreover, one can think of already existing
specifications of complicated finite element problems that via UFC can
be combined with other libraries than the specifications were
originally written for. We have tried to make minimal assumptions to
make this possible.

We believe that UFC itself and the ideas behind it constitute an
important step towards greater flexibility in finite element
software. By code generation via tools like FFC and SFC, this
flexibility may be retained also in combination with very high
performance.

\bibliographystyle{chicago}
\bibliography{bibliography}

\end{document}